\documentclass[leqno,draft]{article}

\def\dist{\mathop{\rm dist}}

\newtheorem{theorem}{Theorem}
\newtheorem{lemma}[theorem]{Lemma}
\newtheorem{proposition}[theorem]{Proposition}
\newtheorem{definition}[theorem]{Definition}
\newtheorem{corollary}[theorem]{Corollary}

\newcommand{\begintheorem}{\addtocounter{equation}{1}\begin{theorem}}
\newcommand{\beginlemma}{\addtocounter{equation}{1}\begin{lemma}}
\newcommand{\beginproposition}{\addtocounter{equation}{1}\begin{proposition}}
\newcommand{\begindefinition}{\addtocounter{equation}{1}\begin{definition}}
\newcommand{\begincorollary}{\addtocounter{equation}{1}\begin{corollary}}

\begin{document}

\title{Some remarks about solenoids, 3}

\author{Stephen Semmes \\
        Rice University}

\date{}

\maketitle

\begin{abstract}
A basic class of constructions is considered, in connection with
bilipschitz mappings in particular.
\end{abstract}

\tableofcontents

\section{A basic situation}
\label{a basic situation}
\setcounter{equation}{0}

        Let $X$ be a (nonempty) Hausdorff topological space, and suppose
that $\phi$ is a homeomorphism from $X$ onto itself.  Thus $X \times
[0, 1]$ is also a Hausdorff space with respect to the product
topology, using the standard topology on the unit interval $[0, 1]$.
Let $\sim_1$ be the equivalence relation on $X \times [0, 1]$ in which
every element of $X \times [0, 1]$ is equivalent to itself, and otherwise
\begin{equation}
\label{(x, 0) sim_1 (phi(x), 1)}
        (x, 0) \sim_1 (\phi(x), 1)
\end{equation}
for every $x \in X$.  This leads to a quotient space
\begin{equation}
\label{Y_1 = (X times [0, 1]) / sim_1}
        Y_1 = (X \times [0, 1]) / \sim_1,
\end{equation}
where the two ends $X \times \{0\}$ and $X \times \{1\}$ of $X \times
[0, 1]$ are glued together using $\phi$.  Let $q_1$ be the
corresponding quotient mapping from $X \times [0, 1]$ onto $Y_1$, so that
\begin{equation}
\label{q_1((x, 0)) = q_1((phi(x), 1))}
        q_1((x, 0)) = q_1((\phi(x), 1))
\end{equation}
for every $x \in X$, and otherwise $q_1$ is one-to-one.  The quotient
topology on $Y_1$ is defined as usual by saying that $U \subseteq Y_1$
is an open set in $Y_1$ if and only if $q_1^{-1}(U)$ is an open set in
$X \times [0, 1]$.  In particular, $q_1$ is automatically continuous
with respect to the quotient topology on $Y_1$, and it is easy to see
that $Y_1$ is also a Hausdorff space under these conditions.  If $X$
is compact, then $X \times {\bf R}$ is compact too, and hence $Y_1$ is
compact with respect to the quotient topology.

        Let $\sim$ be the equivalence relation on $[0, 1]$ in which
every element of $[0, 1]$ is equivalent to itself, and $0$ is equivalent
to $1$.  Thus the quotient topological space
\begin{equation}
\label{[0, 1] / sim}
        [0, 1] / \sim
\end{equation}
is obtained by gluing the ends of $[0, 1]$ together, and is
homeomorphic to the unit circle ${\bf S}^1$ with the standard
topology.  The obvious coordinate projection from $X \times [0, 1]$
onto $[0, 1]$ leads to a continuous mapping from $Y_1$ onto (\ref{[0, 1]
  / sim}), whose fibers are homeomorphic to $X$.  If $\phi$ is the
identity mapping on $X$, then $Y_1$ is homeomorphic to the product of
$X$ and (\ref{[0, 1] / sim}) in a simple way.

        Alternatively, let $\Phi$ be the mapping from $X \times {\bf R}$ 
into itself defined by
\begin{equation}
\label{Phi((x, t)) = (phi(x), t + 1)}
        \Phi((x, t)) = (\phi(x), t + 1)
\end{equation}
for every $x \in X$ and $t \in {\bf R}$.  Thus $\Phi$ is a homeomorphism
from $X \times {\bf R}$ onto itself.  If $n$ is a positive integer, then
\begin{equation}
        \Phi^n((x, t)) = (\phi^n(x), t + n)
\end{equation}
for every $x \in X$ and $t \in {\bf R}$, where $\phi^n$ and $\Phi^n$
are the $n$-fold compositions of these mappings on the correspondng spaces.
This also works for $n = 0$, where the $n$-fold composition is interpreted
as being the identity mapping on the appropriate space, and when $n$ is a
negative integer, for which the $n$-fold composition is considered to be the
$(-n)$-fold composition of the inverse mapping.

        The collection of mappings $\Phi^n$ with $n \in {\bf Z}$ is a
group of homeomorphisms on $X \times {\bf R}$.  This leads to an
equivalence relation $\sim_2$ on $X \times {\bf R}$, where
\begin{equation}
\label{(x, t) sim_2 (x', t')}
        (x, t) \sim_2 (x', t')
\end{equation}
for some $x, x' \in X$ and $t, t' \in {\bf R}$ if and only if there is an
integer $n$ such that
\begin{equation}
\label{Phi^n((x, t)) = (x', t')}
        \Phi^n((x, t)) = (x', t').
\end{equation}
Let $q_2$ be the quotient mapping from $X \times {\bf R}$ onto the
quotient space
\begin{equation}
\label{Y_2 = (X times {bf R}) / sim_2}
        Y_2 = (X \times {\bf R}) / \sim_2.
\end{equation}
As before, the quotient topology on $Y_2$ is defined by saying that $U
\subseteq Y_2$ is an open set if and only if $q_2^{-1}(U)$ is an open
set in $X \times {\bf R}$, so that the quotient mapping $q_2$ is 
automatically continuous.

        Let us consider the restriction of $q_2$ to $X \times [0, 1]
\subseteq X \times {\bf R}$.  By construction, if $x, x' \in X$ and 
$t, t' \in [0, 1]$, then
\begin{equation}
\label{q_1((x, t)) = q_1((x', t'))}
        q_1((x, t)) = q_1((x', t'))
\end{equation}
in $Y_1$ if and only if
\begin{equation}
\label{q_2((x, t)) = q_2((x', t'))}
        q_2((x, t)) = q_2((x', t'))
\end{equation}
in $Y_2$.  This leads to a mapping from $Y_1$ into $Y_2$, which is
easily seen to be a homeomorphism from $Y_1$ onto $Y_2$.  An advantage
of $Y_2$ is that $q_2$ is a local homeomorphism from $X \times {\bf
  R}$ onto $Y_2$.  Although $q_1$ is a local homeomorphism around
points $(x, t) \in X \times (0, 1)$, this does not work when $t = 0$
or $1$.

        Of course, the real line ${\bf R}$ is a commutative topological 
group with respect to addition, which contains ${\bf Z}$ as a discrete
subgroup.  The quotient ${\bf R} / {\bf Z}$ is also a commutative
topological group with respect to the quotient topology and group
operation, which is isomorphic as a topological group to the
multiplicative group of complex numbers with modulus $1$.  The obvious
coordinate projection from $X \times {\bf R}$ onto ${\bf R}$ leads to
a continuous mapping from $Y_2$ onto ${\bf R} / {\bf Z}$.  This mapping
corresponds exactly to the continuous mapping from $Y_1$ onto 
(\ref{[0, 1] / sim}) discussed earlier, using the identification between
$Y_1$ and $Y_2$ described in the previous paragraph.  This also uses
the analogous identification between (\ref{[0, 1] / sim}) and
${\bf R} / {\bf Z}$.

        If $r \in R$, then
\begin{equation}
\label{(x, t) mapsto (x, t + r)}
        (x, t) \mapsto (x, t + r)
\end{equation}
defines a homeomorphism from $X \times {\bf R}$ onto itself that
preserves the equivalence relation $\sim_2$.  This leads to a 
homeomorphism $A_r$ from $Y_2$ onto itself, where
\begin{equation}
\label{A_r(q_2((x, t))) = q_2((x, t + r))}
        A_r(q_2((x, t))) = q_2((x, t + r))
\end{equation}
for every $x \in X$ and $t \in {\bf R}$.  This is actually a group of
homeomorphisms from $Y_2$ onto itself, in the sense that
\begin{equation}
\label{A_r circ A_r' = A_{r + r'}}
        A_r \circ A_r' = A_{r + r'}
\end{equation}
for each $r, r' \in {\bf R}$, because of the analogous property of
(\ref{(x, t) mapsto (x, t + r)}) on $X \times {\bf R}$.  Note that
\begin{equation}
\label{A_n(q_2((phi^n(x), t))) = q_2(x, t))}
        A_n(q_2((\phi^n(x), t))) = q_2((x, t))
\end{equation}
for each $x \in X$, $t \in {\bf R}$, and $n \in {\bf Z}$.

        It is easy to see that
\begin{equation}
\label{psi_t(x) = q_2((x, t))}
        \psi_t(x) = q_2((x, t))
\end{equation}
defines a homeomorphism from $X$ onto $q_2(X \times \{t\})$ for each
$t \in {\bf R}$, where $q_2(X \times \{t\})$ is equipped with the
topology induced by the one on $Y_2$.  The sets $q_2(X \times \{t\})$
with $t \in {\bf R}$ are the fibers of the natural projection from
$Y_2$ onto ${\bf R} / {\bf Z}$, which satisfy the periodicity condition
\begin{equation}
\label{q_2(X times {t + 1}) = q_2(X times {t})}
        q_2(X \times \{t + 1\}) = q_2(X \times \{t\})
\end{equation}
for each $t \in {\bf R}$, by construction.  More precisely,
\begin{equation}
\label{psi_{t + 1}(phi(x)) = q_2((phi(x), t + 1)) = q_2(x, t) = psi_t(x)}
        \psi_{t + 1}(\phi(x)) = q_2((\phi(x), t + 1)) = q_2(x, t) = \psi_t(x)
\end{equation}
for every $x \in X$ and $t \in {\bf R}$, which implies that $\psi_{t +
  1}(X) = \psi_t(X)$.  Similarly,
\begin{equation}
\label{A_r(psi_t(x)) = psi_{r + t}(x)}
        A_r(\psi_t(x)) = \psi_{r + t}(x)
\end{equation}
for every $x \in X$ and $r, t \in {\bf R}$, and $A_r$ maps $q_2(X
\times \{t\})$ onto $q_2(X \times \{r + t\})$ for each $r, t \in {\bf
  R}$.

\section{Connectedness}
\label{connectedness}
\setcounter{equation}{0}

        Let us continue with the notation and hypotheses in the previous
section.  Let $x \in X$ be given, and consider
\begin{equation}
\label{q_2({x} times {bf R}) = {A_r((x, 0)) : r in {bf R}}}
        q_2(\{x\} \times {\bf R}) = \{A_r((x, 0)) : r \in {\bf R}\}.
\end{equation}
If $\phi^n(x) \ne x$ for every positive integer $n$, then it is easy
to see that the restriction of $q_2$ to $\{x\} \times {\bf R}$ is a
one-to-one mapping into $Y_2$.  Otherwise, if $\phi^n(x) = x$ for some
$x \in {\bf Z}_+$, then $q_2((x, t))$ is periodic in $t$, with period
$n$.  Note that (\ref{q_2({x} times {bf R}) = {A_r((x, 0)) : r in {bf
      R}}}) is a connected set in $Y_2$ for each $x \in X$, because
the real line is connected.  If $X$ is totally disconnected, then the
subsets of $X \times {\bf R}$ of the form $\{x\} \times {\bf R}$ for
some $x \in X$ are the pathwise-connected components of $X \times {\bf
  R}$.  In this case, the subsets of $Y_2$ of the form (\ref{q_2({x}
  times {bf R}) = {A_r((x, 0)) : r in {bf R}}}) for some $x \in X$
are the pathwise-connected components of $Y_2$.  Of course, if $X$
is connected, then $X \times {\bf R}$ is connected, and hence
$Y_2 = q_2(X \times {\bf R})$ is connected too.

        If $E \subseteq Y_2$ is both open and closed, then it follows that 
for each $x \in X$, (\ref{q_2({x} times {bf R}) = {A_r((x, 0)) : r in
    {bf R}}}) is either contained in $E$ or in $Y_2 \setminus E$.
Equivalently, this means that
\begin{equation}
\label{A_r(E) = E}
        A_r(E) = E
\end{equation}
for each $r \in {\bf R}$, so that $E$ is invariant under the flow on $Y_2$
defined by $A_r$.  Put
\begin{equation}
\label{E_0 = {x in X : q_2({x} times {bf R}) subseteq E}}
        E_0 = \{x \in X : q_2(\{x\} \times {\bf R}) \subseteq E\},
\end{equation}
and observe that $\phi(E_0) = E_0$, since
\begin{equation}
\label{q_2({phi(x)} times {bf R}) = q_2({x} times {bf R})}
        q_2(\{\phi(x)\} \times {\bf R}) = q_2(\{x\} \times {\bf R})
\end{equation}
for each $x \in X$.  Alternatively,
\begin{equation}
        E_0 \times {\bf R} = q_2^{-1}(E),
\end{equation}
which is automatically invariant under $\Phi$, and
\begin{equation}
\label{E_0 = {x in X : psi_t(x) in E}}
        E_0 = \{x \in X : \psi_t(x) \in E\}
\end{equation}
for each $t \in {\bf R}$.  This implies that that $E_0$ is both open
and closed in $X$.  Note too that $E_0 \ne \emptyset$ when $E \ne
\emptyset$, and that $E_0 \ne X$ when $E \ne Y_2$.  It follows that if
$Y_2$ is not connected, then there is an open and closed set $E_0
\subseteq X$ such that $E_0 \ne \emptyset, X$ and $\phi(E_0) = E_0$.

        Conversely, suppose that $E_0 \subseteq X$ is both open and
closed in $X$, and that $\phi(E_0) = E_0$.  This implies that
\begin{equation}
\label{Phi(E_0 times {bf R}) = E_0 times {bf R}}
        \Phi(E_0 \times {\bf R}) = E_0 \times {\bf R},
\end{equation}
and we put
\begin{equation}
\label{E = q_2(E_0 times {bf R})}
        E = q_2(E_0 \times {\bf R}),
\end{equation}
which is automatically invariant under $A_r$ for each $r \in {\bf R}$. 
Observe that
\begin{equation}
\label{Y_2 setminus E = q_2((X setminus E_0) times {bf R})}
        Y_2 \setminus E = q_2((X \setminus E_0) \times {\bf R}),
\end{equation}
because $E_0 \times {\bf R}$ is invariant under $\Phi$, and hence that
$E$ is both open and closed in $Y_2$.  If $E_0 \ne \emptyset, X$, then
$E \ne \emptyset, Y_2$, and thus $Y_2$ is not connected.  This shows
that $Y_2$ is connected if and only if there is no set $E_0 \subseteq X$
such that $E_0$ is both open and closed in $X$, $E_0 \ne \emptyset, X$,
and $\phi(E_0) = E_0$.

        If $x_0 \in X$ and the orbit
\begin{equation}
\label{{phi^n(x_0) : n in {bf Z}}}
        \{\phi^n(x_0) : n \in {\bf Z}\}
\end{equation}
of $x_0$ under $\phi$ is dense in $X$, then one can check that
$q_2(\{x_0\} \times {\bf R})$ is dense in $Y_2$.  This implies that
$Y_2$ is connected, since the closure of a connected set is connected.
Alternatively, if $E_0 \subseteq X$ satisfies $\phi(E_0) = E_0$, then
the orbit of every element of $X$ under $\phi$ is either contained in
$E_0$ or in $X \backslash E_0$.  If $E_0$ is also both open and closed
in $X$, then $E_0$ and $X \backslash E_0$ are both closed sets in $X$,
and hence the closure of the orbit of every element of $X$ under
$\phi$ is contained in $E_0$ of $X \backslash E_0$.  If additionally
$E_0 \ne \emptyset, X$, so that $E_0$ and $X \backslash E_0$ are both
proper subsets of $X$, then it follows that the closure of the orbit of
any element of $X$ under $\phi$ is proper subset of $X$ as well.

\section{Topological groups}
\label{topological groups}
\setcounter{equation}{0}

        Let $G$ be a topological group, and let $h$ be an element of
$G$.  Thus
\begin{equation}
\label{phi(x) = x h}
        \phi(x) = x \, h
\end{equation}
defines a homeomorphism from $G$ onto itself, and
\begin{equation}
\label{phi^n(x) = x h^n}
        \phi^n(x) = x \, h^n
\end{equation}
for each $n \in {\bf Z}$.  Note that $G \times {\bf R}$ is also a
topological group, where the group operations are defined
coordinatewise, and using the product topology.  Let $H$ be the
subgroup of $G \times {\bf R}$ consisting of $(h^n, n)$ for each
integer $n$, which is a discrete subgroup of $G \times {\bf R}$.  Thus
the quotient space $(G \times {\bf R}) / H$ of left cosets of $H$ in
$G \times {\bf R}$ can be defined in the usual way, with the quotient
topology on $(G \times {\bf R}) / H$ associated to the product
topology on $G \times {\bf R}$.  The quotient space $(G \times {\bf
  R}) / H$ corresponds exactly to the space $Y_2$ in Section \ref{a
  basic situation}, and the natural quotient mapping from $G \times
{\bf R}$ onto $(G \times {\bf R}) / H$ corresponds to the mapping
$q_2$ in Section \ref{a basic situation}.  If the subgroup of $G$
generated by $h$ is normal, then $H$ is a normal subgroup in $G \times
{\bf R}$, and $(G \times {\bf R}) / H$ is a topological group as well.
Otherwise, $G \times {\bf R}$ acts on the quotient space $(G \times
{\bf R}) / H$ by left translations.  If $h$ is the identity element in
$G$, then $(G \times {\bf R}) / H$ reduces to $G \times ({\bf R} /
{\bf Z})$.

        Of course, the subgroup of $G$ generated by $h$ is abelian,
and hence its closure in $G$ is abelian.  In particular, if the
subgroup of $G$ generated by $h$ is dense in $G$, then $G$ is abelian.
This would also imply that $(G \times {\bf R}) / H$ is connected, as
in the previous section.  If $G = {\bf Z}$ as a discrete group with
respect to addition and $h = 1$, then it is easy to see that $(G
\times {\bf R}) / H$ is isomorphic as a topological group to ${\bf
  R}$.  Alternatively, let $p$ be a prime number, and let ${\bf Z}_p$
be the group of $p$-adic integers.  This is a compact totally disconnected
commutative topological group with respect to addition, which contains
${\bf Z}$ as a dense subgroup.  If we take $h = 1$ as an element of
${\bf Z}_p$, then the correspondng quotient $(G \times {\bf R}) / H$
is a compact commutative topological group which is connected but not
locally connected.

        If there is a countable local base for the topology of $G$
at the identity element, then a famous theorem states that there
is a metric on $G$ that determines the same topology and which is
invariant under right translations.  We shall look at isometric
mappings more broadly in the next section.

\section{Isometries}
\label{isometries}
\setcounter{equation}{0}

        Let us return now to the setting of Section \ref{a basic situation}.
Suppose in addition that the topology on $X$ is determined by a metric 
$d(x, y)$, and that $\phi$ is an isometric mapping from $X$ onto itself,
so that
\begin{equation}
\label{d(phi(x), phi(y)) = d(x, y)}
        d(\phi(x), \phi(y)) = d(x, y)
\end{equation}
for every $x, x' \in X$.  Put
\begin{equation}
\label{rho((x, r), (y, t)) = max(d(x, y), |r - t|)}
        \rho((x, r), (y, t)) = \max(d(x, y), |r - t|)
\end{equation}
for each $x, y \in X$ and $r, t \in {\bf R}$, which defines a metric
on $X \times {\bf R}$ for which the corresponding topology is the
product topology.  Thus
\begin{equation}
\label{rho(Phi((x, r)), Phi((y, t))) = ... = rho((x, r), (y, t))}
  \quad \rho(\Phi((x, r)), \Phi((y, t))) 
          = \rho((\phi(x), r + 1), (\phi(y), t + 1)) = \rho((x, r), (y, t))
\end{equation}
for every $x, y \in X$ and $r, t \in {\bf R}$, where $\Phi$ is the
mapping from $X \times {\bf R}$ onto itself defined in Section \ref{a
  basic situation}.

        The corresponding quotient metric on $Y_2$ is defined by
\begin{eqnarray}
\label{D(q_2((x, r)), q_2((y, t))) = ...}
\lefteqn{D(q_2((x, r)), q_2((y, t)))} \\
     & = & \inf \{\rho((x', r'), (y', t')) : x', y' \in X, \ 
                                           r', t' \in {\bf R}, \nonumber \\
     & & \qquad q_2((x', r')) = q_2((x, r)), \ q_2((y', t')) = q_2((y, t))\} 
                                                                \nonumber
\end{eqnarray}
for each $x, y \in X$ and $r, t \in {\bf R}$.  Equivalently,
\begin{eqnarray}
\label{D(q_2((x, r)), q_2((y, t))) = ..., 2}
\lefteqn{D(q_2((x, r)), q_2((y, t)))} \\
 & = & \inf\{\rho((x', r'), (y, t)) : x' \in X, \ r' \in {\bf R}, \ 
                                     q_2((x', r')) = q_2((x, r))\} \nonumber \\
 & = & \inf\{\rho((x, r), (y', t')) : y' \in X, \ t' \in {\bf R}, \ 
                                     q_2((y', t')) = q_2((y, t))\}, \nonumber
\end{eqnarray}
because $\Phi$ is an isometry on $X \times {\bf R}$ with respect to
$\rho(\cdot, \cdot)$.  If $x, x', y, z, z' \in X$ and $r, r', t, u, u'
\in {\bf R}$ satisfy $q_2((x', r')) = q_2((x, r))$ and $q_2((z', u')) =
q_2((z, u))$, then
\begin{eqnarray}
\label{D(q_2((x, r)), q_2((z, u))) le rho((x', r'), (z', u')) le ...}
        D(q_2((x, r)), q_2((z, u))) & \le & \rho((x', r'), (z', u')) \\
           & \le & \rho((x', r'), (y, t)) + \rho((y, t), (z', u')), \nonumber
\end{eqnarray}
by the triangle inequality for $\rho(\cdot, \cdot)$.  Taking the infimum
over $(x', r')$ and $(z', u')$, we get that
\begin{eqnarray}
\label{D(q_2((x, r)), q_2((z, u))) le ...}
\lefteqn{D(q_2((x, r)), q_2((z, u)))} \\
 & \le & D(q_2((x, r)), q_2((y, t))) + D(q_2((y, t)), q_2((z, u))). \nonumber
\end{eqnarray}
Thus $D(\cdot, \cdot)$ satisfies the triangle inequality on $Y_2$, and
it is easy to see that $D(\cdot, \cdot)$ is a metric on $Y_2$ that
defines the same topology on $Y_2$ as before.

        By construction,
\begin{equation}
\label{D(q_2((x, r)), q_2((y, t))) le rho((x, r), (y, t))}
        D(q_2((x, r)), q_2((y, t))) \le \rho((x, r), (y, t))
\end{equation}
for every $x, y \in X$ and $r, t \in {\bf R}$.  Suppose that $r, t \in
{\bf R}$ satisfy
\begin{equation}
\label{|r - t| le 1/2}
        |r - t| \le 1/2,
\end{equation}
so that
\begin{equation}
\label{|r - t'| ge 1/2}
        |r - t'| \ge 1/2
\end{equation}
for every $t' \in {\bf R}$ such that $t' - t \in {\bf Z}$ and $t' \ne t$.
This implies that
\begin{equation}
\label{rho((x, r), (y', t')) ge |r - t'| ge 1/2}
        \rho((x, r), (y', t')) \ge |r - t'| \ge 1/2
\end{equation}
for every $x, y, y' \in X$ and $t' \in {\bf R}$ such that $q_2((y', t'))
= q_2((y, t))$ and $(y', t') \ne (y, t)$, so that
\begin{equation}
\label{D(q_2((x, r)), q_2((y, t))) ge min(rho((x, r), (y, t)), 1/2)}
 D(q_2((x, r)), q_2((y, t))) \ge \min\big(\rho((x, r), (y, t)), 1/2\big)
\end{equation}
by (\ref{D(q_2((x, r)), q_2((y, t))) = ..., 2}).  In particular, if
\begin{equation}
\label{d(x, y) le 1/2}
        d(x, y) \le 1/2,
\end{equation}
then $\rho((x, r), (y, t)) \le 1/2$, and hence
\begin{equation}
\label{D(q_2((x, r)), q_2((y, t))) = rho((x, r), (y, t))}
        D(q_2((x, r)), q_2((y, t))) = \rho((x, r), (y, t)),
\end{equation}
by (\ref{D(q_2((x, r)), q_2((y, t))) le rho((x, r), (y, t))}) and
(\ref{D(q_2((x, r)), q_2((y, t))) ge min(rho((x, r), (y, t)), 1/2)}).
Similarly, if
\begin{equation}
\label{d(x, y) le k}
        d(x, y) \le k
\end{equation}
for some $k \ge 1/2$, then $\rho((x, r), (y, t)) \le k$, and we get that
\begin{equation}
\label{rho((x, r), (y, t)) le 2 k D(q_2((x, r)), q_2((y, t)))}
        \rho((x, r), (y, t)) \le 2 \, k \, D(q_2((x, r)), q_2((y, t))).
\end{equation}

        If $X$ is bounded with respect to $d(x, y)$, then (\ref{d(x, y) le k})
holds for some $k \ge 1/2$ and every $x, y \in X$.  This implies that
(\ref{rho((x, r), (y, t)) le 2 k D(q_2((x, r)), q_2((y, t)))}) holds
for every $x, y \in X$ and $r, t \in {\bf R}$ that satisfy (\ref{|r -
  t| le 1/2}).  Otherwise, for any positive real number $k$,
\begin{equation}
\label{d_1(x, y) = min(d(x, y), k)}
        d_1(x, y) = \min(d(x, y), k)
\end{equation}
defines a metric on $X$ which is topologically equivalent to $d(x,
y)$.  Of course, if $\phi$ is an isometry on $X$ with respect to $d(x,
y)$, then $\phi$ is an isometry on $X$ with respect to $d_1(x, y)$ as
well.

        Suppose now that $\phi$ is not necessarily an isometry on $X$ with 
respect to $d(x, y)$, but that the collection of iterates $\phi^n$ with 
$n \in {\bf Z}$ is equicontinuous at every point in $X$ with respect to 
$d(x, y)$.  This means that for each $x \in X$ and $\epsilon > 0$ there is 
a $\delta(x, \epsilon) > 0$ such that
\begin{equation}
\label{d(phi^n(x), phi^n(y)) le epsilon}
        d(\phi^n(x), \phi^n(y)) \le \epsilon
\end{equation}
for every $n \in {\bf Z}$ and $y \in X$ such that $d(x, y) < \delta(x,
\epsilon)$.  We may as well ask that $X$ be bounded with respect to
$d(x, y)$ too, since otherwise we can replace $d(x, y)$ with
(\ref{d_1(x, y) = min(d(x, y), k)}) for some $k > 0$, and still have
the same equicontinuity condition.  If we put
\begin{equation}
\label{widetilde{d}(x, y) = sup_{n in {bf Z}} d(phi^n(x), phi^n(y))}
        \widetilde{d}(x, y) = \sup_{n \in {\bf Z}} d(\phi^n(x), \phi^n(y)),
\end{equation}
then $\widetilde{d}(x, y)$ is a metric on $X$, and
\begin{equation}
\label{d(x, y) le widetilde{d}(x, y)}
        d(x, y) \le \widetilde{d}(x, y)
\end{equation}
for every $x, y \in X$, since we can take $n = 0$ in
(\ref{widetilde{d}(x, y) = sup_{n in {bf Z}} d(phi^n(x), phi^n(y))}).
We also have that
\begin{equation}
\label{widetilde{d}(x, y) le epsilon}
        \widetilde{d}(x, y) \le \epsilon
\end{equation}
for every $x, y \in X$ such that $d(x, y) < \delta(x, \epsilon)$, by
(\ref{d(phi^n(x), phi^n(y)) le epsilon}), and hence that
$\widetilde{d}(x, y)$ and $d(x, y)$ determine the same topology on
$X$.  By construction,
\begin{equation}
\label{widetilde{d}(phi(x), phi(y)) = widetilde{d}(x, y)}
        \widetilde{d}(\phi(x), \phi(y)) = \widetilde{d}(x, y)
\end{equation}
for every $x, y \in X$, so that $\phi$ is an isometry on $X$ with
respect to $\widetilde{d}$.  This is a bit nicer when the collection
of iterates $\phi^n$ with $n \in {\bf Z}$ is uniformly equicontinuous
on $X$, in the sense that one can take $\delta(x, \epsilon) =
\delta(\epsilon)$ independent of $x \in X$ for each $\epsilon > 0$.
In this case, the identity mapping on $X$ is uniformly continuous as a
mapping from $X$ equipped with $d(x, y)$ onto $X$ equipped with
$\widetilde{d}(x, y)$.

        Suppose that $X$ is bounded with respect to $d(x, y)$, and
let $C(X, X)$ be the space of continuous mappings from $X$ onto itself.
Thus the supremum metric
\begin{equation}
\label{sigma(f, g) = sup_{x in X} d(f(x), g(x))}
        \sigma(f, g) = \sup_{x \in X} d(f(x), g(x))
\end{equation}
is defined for each $f, g \in C(X, X)$, and it is easy to see that the
group $\mathcal{I}(X)$ of isometric mappings from $X$ onto itself is a
topological group with respect to the topology determined by the
restriction of $\sigma(f, g)$ to $\mathcal{I}(X)$.  If $X$ is compact,
then one can check that $\mathcal{I}(X)$ is compact with respect to
the supremum metric, using the Arzela--Ascoli theorem.

        It is easy to see that the distance
\begin{equation}
\label{dist(a, {bf Z}) = min_{n in {bf Z}} |a - n|}
        \dist(a, {\bf Z}) = \min_{n \in {\bf Z}} |a - n|
\end{equation}
from $a \in {\bf R}$ to ${\bf Z}$ satisfies
\begin{equation}
\label{dist(a + b, {bf Z}) le dist(a, {bf Z}) + dist(b, {bf Z})}
        \dist(a + b, {\bf Z}) \le \dist(a, {\bf Z}) + \dist(b, {\bf Z})
\end{equation}
for every $a, b \in {\bf R}$.  Note that $\dist(r - t, {\bf Z})$ is
the same as the distance between the images of $r, t \in {\bf R}$ in
${\bf R} / {\bf Z}$ under the natural quotient mapping from ${\bf R}$
onto ${\bf R} / {\bf Z}$, with respect to the quotient metric on ${\bf
  R} / {\bf Z}$ associated to the standard metric on ${\bf R}$.  Of
course,
\begin{equation}
\label{rho((x, r), (y, t)) ge |r - t|}
        \rho((x, r), (y, t)) \ge |r - t|
\end{equation}
for every $x, y \in X$ and $r, t \in {\bf R}$, by construction.  
It follows that
\begin{equation}
\label{D(q_2((x, r)), q_2((y, t))) ge dist(r - t, {bf Z})}
        D(q_2((x, r)), q_2((y, t))) \ge \dist(r - t, {\bf Z})
\end{equation}
for every $x, y \in X$ and $r, t \in {\bf R}$, by the definition of
$D(q_2((x, r)), q_2((y, t)))$.

\section{Bilipschitz mappings}
\label{bilipschitz mappings}
\setcounter{equation}{0}

        Let us go back to the setting of Section \ref{a basic situation},
and suppose again that the topology on $X$ is determined by a metric
$d(x, y)$.  Instead of asking that $\phi$ be an isometry on $X$, let us
suppose that $\phi$ is bilipschitz, so that
\begin{equation}
\label{C^{-1} d(x, y) le d(phi(x), phi(y)) le C d(x, y)}
        C^{-1} \, d(x, y) \le d(\phi(x), \phi(y)) \le C \, d(x, y)
\end{equation}
for some $C \ge 1$ and every $x, y \in X$.  Of course, this implies
that $\phi$ is an isometry on $X$ when $C = 1$.  Otherwise, note that
$\phi^{-1}$ is also bilipschitz with the same constant $C$, and that
$\phi^n$ is bilipschitz with constant $C^{|n|}$ for each $n \in {\bf
  Z}$.  If $\phi^n$ is actually bilipschitz with a constant that does
not depend on $n$ for each $n \in {\bf Z}$, then $\phi$ is an isometry
with respect to the metric $\widetilde{d}(x, y)$ on $X$ defined in the
previous section, and $\widetilde{d}(x, y)$ is bounded by a constant
multiple of $d(x, y)$.

        As in the previous section,
\begin{equation}
\label{rho((x, r), (y, t)) = max(d(x, y), |r - t|), 2}
        \rho((x, r), (y, t)) = \max(d(x, y), |r - t|)
\end{equation}
defines a metric on $X \times {\bf R}$ for which the corresponding
topology is the product topology.  Let $\Phi$ be the mapping on
$X \times {\bf R}$ defined in Section \ref{a basic situation},
so that
\begin{eqnarray}
\label{rho(Phi(x, r), Phi(y, t)) = ... = max(d(phi(x), phi(y)), |r - t|)}
 \rho(\Phi(x, r), \Phi(y, t)) & = & \rho((\phi(x), r + 1), (\phi(y), t + 1)) \\
                         & = & \max(d(\phi(x), \phi(y)), |r - t|) \nonumber
\end{eqnarray}
for every $x, y \in X$ and $r, t \in {\bf R}$.  Using this, it is easy
to see that $\Phi$ is a bilipschitz mapping on $X \times {\bf R}$ with
constant $C$ with respect to $\rho(\cdot, \cdot)$.  As before, we
would like to define a distance function
\begin{equation}
\label{delta(q_2((x, r)), q_2((y, t)))}
        \delta(q_2((x, r)), q_2((y, t)))
\end{equation}
on $Y_2$ that looks locally like (\ref{rho((x, r), (y, t)) = max(d(x,
  y), |r - t|), 2}), at least when $|r|$ and $|t|$ are not too large.

        Let $x, y \in X$ and $r, t \in {\bf R}$ be given, and suppose that
$x', y' \in X$ and $r', t' \in {\bf R}$ satisfy
\begin{equation}
\label{q_2((x, r)) = q_2((x', r')), q_2((y, t)) = q_2((y', t'))}
        q_2((x, r)) = q_2((x', r')), \quad q_2((y, t)) = q_2((y', t')).
\end{equation}
This implies that
\begin{equation}
\label{r equiv r' and t equiv t' modulo {bf Z}}
        r \equiv r' \hbox{ and } t \equiv t' \hbox{ modulo } {\bf Z},
\end{equation}
and hence
\begin{equation}
\label{r - t equiv r' - t' modulo {bf Z}}
        r - t \equiv r' - t' \hbox{ modulo } {\bf Z}.
\end{equation}
Note that $x'$, $y'$ are uniquely determined by (\ref{q_2((x, r)) =
  q_2((x', r')), q_2((y, t)) = q_2((y', t'))}) and $r'$, $t'$.  If we
restrict our attention to $r'$, $t'$ in a bounded set, then there are
only finitely many possibilities for them, and thus for $x'$, $y'$.

        We can always choose $r', t' \in {\bf R}$ so that 
(\ref{q_2((x, r)) = q_2((x', r')), q_2((y, t)) = q_2((y', t'))}) holds and
\begin{equation}
\label{|r' - t'| le 1/2}
        |r' - t'| \le 1/2,
\end{equation}
by adding suitable integers to $r'$ or $t'$.  One can also get
\begin{equation}
\label{|frac{r' + t'}{2}| le frac{1}{2}}
        \biggl|\frac{r' + t'}{2}\biggr| \le \frac{1}{2},
\end{equation}
by adding a suitable integer to both $r'$ and $t'$, which does not
affect (\ref{|r' - t'| le 1/2}).  Under these conditions,
\begin{equation}
\label{|r'|, |t'| le 3/4}
        |r'|, |t'| \le 3/4,
\end{equation}
because the distance from $r'$ or $t'$ to $(r' + t')/2$ is equal to
$|r' - t'|/2$.

        Put
\begin{eqnarray}
\label{delta(q_2((x, r)), q_2((y, t))) = ...}
\lefteqn{\delta(q_2((x, r)), q_2((y, t)))}           \\
 & = & \min\{\rho((x', r'), (y', t')) : x', y' \in X, \ r', t' \in {\bf R}
                                                         \nonumber \\
  & & \qquad\qquad\qquad\quad \hbox{satisfy } 
        (\ref{q_2((x, r)) = q_2((x', r')), q_2((y, t)) = q_2((y', t'))}), 
 (\ref{|r' - t'| le 1/2}), \hbox{ and } (\ref{|r'|, |t'| le 3/4})\}. \nonumber
\end{eqnarray}
As in the previous paragraphs, every pair of points in $Y_2$ can be
represented in this way, and there are only finitely many such
representations.  Thus the minimum in (\ref{delta(q_2((x, r)), q_2((y,
  t))) = ...}) makes sense, and is a nonnegative real number.  If
\begin{equation}
        q_2((x, r)) = q_2((y, t)),
\end{equation}
then we can choose $x', y' \in X$ and $r', t' \in {\bf R}$ such that
$|r'| \le 1/2 \le 1$ and $r' = t'$, which implies that
(\ref{delta(q_2((x, r)), q_2((y, t))) = ...}) is equal to $0$.
Otherwise, if $q_2((x, r)) \ne q_2((y, t))$, then (\ref{delta(q_2((x,
  r)), q_2((y, t))) = ...}) is the minimum of finitely many positive
real numbers, and hence is positive too.  Clearly (\ref{delta(q_2((x,
  r)), q_2((y, t))) = ...}) is symmetric in $q_2((x, r))$ and $q_2((y,
t))$.  However, (\ref{delta(q_2((x, r)), q_2((y, t))) = ...}) does not
normally satisfy the triangle inequality, and we shall come back to
that soon.

        Suppose that
\begin{equation}
\label{|r|, |t| le 3/4 and |r - t| le 1/2}
        |r|, |t| \le 3/4 \quad \hbox{and} \quad |r - t| \le 1/2,
\end{equation}
so that
\begin{equation}
\label{delta(q_2((x, r)), q_2((y, t))) le rho((x, r), (y, t))}
        \delta(q_2((x, r)), q_2((y, t))) \le \rho((x, r), (y, t)),
\end{equation}
because $x$, $y$, $r$, and $t$ are admissible competitors for the
minimum in (\ref{delta(q_2((x, r)), q_2((y, t))) = ...}).  If
\begin{equation}
\label{|r - t| < 1/2}
        |r - t| < 1/2,
\end{equation}
then we also have that
\begin{equation}
\label{delta(q_2((x, r)), q_2((y, t))) ge C^{-1} rho((x, r), (y, t))}
        \delta(q_2((x, r)), q_2((y, t))) \ge C^{-1} \, \rho((x, r), (y, t)),
\end{equation}
where $C$ is as in (\ref{C^{-1} d(x, y) le d(phi(x), phi(y)) le C d(x, y)}).
To see this, suppose that $x', y' \in X$ and $r', t' \in {\bf R}$ satisfy 
(\ref{q_2((x, r)) = q_2((x', r')), q_2((y, t)) = q_2((y', t'))}),
(\ref{|r' - t'| le 1/2}), and (\ref{|r'|, |t'| le 3/4}), and that $r' \ne r$
or $t' \ne t$.  Observe that
\begin{equation}
\label{r' - t' = r - t}
        r' - t' = r - t
\end{equation}
in this situation, because of (\ref{r - t equiv r' - t' modulo {bf Z}}),
(\ref{|r' - t'| le 1/2}), and (\ref{|r - t| < 1/2}).  Moreover,
\begin{equation}
\label{|r' - r|, |t' - t| le 3/2}
        |r' - r|, |t' - t| \le 3/2,
\end{equation}
by (\ref{|r'|, |t'| le 3/4}) and (\ref{|r|, |t| le 3/4 and |r - t| le 1/2}),
which implies that
\begin{equation}
\label{|r' - r|, |t' - t| le 1}
        |r' - r|, |t' - t| \le 1,
\end{equation}
because of (\ref{r equiv r' and t equiv t' modulo {bf Z}}).  Thus
\begin{equation}
\label{r' - r = t' - t = 1 or -1}
        r' - r = t' - t = 1 \hbox{ or } -1
\end{equation}
under these conditions, by (\ref{r equiv r' and t equiv t' modulo {bf
    Z}}), (\ref{r' - t' = r - t}), (\ref{|r' - r|, |t' - t| le 1}),
and the hypothesis that $r' \ne r$ or $t' \ne t$.  This implies that
either $x' = \phi(x)$ and $y' = \phi(y)$, or $x' = \phi^{-1}(x)$ and
$y' = \phi^{-1}(y)$, because of (\ref{q_2((x, r)) = q_2((x', r')),
  q_2((y, t)) = q_2((y', t'))}).  In both cases, we get that
\begin{equation}
\label{d(x', y') ge C^{-1} d(x, y)}
        d(x', y') \ge C^{-1} \, d(x, y),
\end{equation}
by (\ref{C^{-1} d(x, y) le d(phi(x), phi(y)) le C d(x, y)}).  It
follows that
\begin{equation}
\label{rho((x', r'), (y', t')) ge C^{-1} rho((x, r), (y, t))}
        \rho((x', r'), (y', t')) \ge C^{-1} \, \rho((x, r), (y, t)),
\end{equation}
using also (\ref{r' - t' = r - t}).  This shows that
(\ref{delta(q_2((x, r)), q_2((y, t))) ge C^{-1} rho((x, r), (y, t))})
holds when $|r - t| < 1/2$, as desired.

        By construction,
\begin{equation}
\label{delta(q_2((x, r)), q_2((y, t))) ge dist(r - t, {bf Z})}
        \delta(q_2((x, r)), q_2((y, t))) \ge \dist(r - t, {\bf Z})
\end{equation}
for every $x, y \in X$ and $r, t \in {\bf R}$, and hence
\begin{equation}
\label{delta(q_2((x, r)), q_2((y, t))) ge 1/2}
        \delta(q_2((x, r)), q_2((y, t))) \ge 1/2
\end{equation}
when $|r - t| = 1/2$.  Combining this with (\ref{delta(q_2((x, r)),
  q_2((y, t))) ge C^{-1} rho((x, r), (y, t))}), we get that
\begin{equation}
\label{delta(q_2((x, r)),q_2((y, t))) ge min(C^{-1} rho((x, r), (y, t)), 1/2)}
        \delta(q_2((x, r)),q_2((y, t)))
                  \ge \min\big(C^{-1} \, \rho((x, r), (y, t)), 1/2\big)
\end{equation}
when $r$, $t$ satisfy (\ref{|r|, |t| le 3/4 and |r - t| le 1/2}).

        As mentioned earlier, $\delta(\cdot, \cdot)$ does not normally
satisfy the triangle inequality.  To fix this, let $q_2((x, r))$
and $q_2((y, t))$ be any two elements of $Y_2$, and consider all
finite sequences
\begin{equation}
\label{q_2((x_1, r_1)), ldots, q_2((x_{n + 1}, r_{n + 1}))}
        q_2((x_1, r_1)), \ldots, q_2((x_{n + 1}, r_{n + 1}))
\end{equation}
of elements of $Y_2$ connecting $q_2((x, r))$ to $q_2((y, t))$, in the
sense that
\begin{equation}
\label{q_2((x_1, r_1)) = q_2((x, r)), q_2(x_{n + 1}, r_{n + 1})) = q_2((y, t))}
        q_2((x_1, r_1)) = q_2((x, r)) \quad \hbox{and} \quad 
                           q_2(x_{n + 1}, r_{n + 1})) = q_2((y, t)).
\end{equation}
Put
\begin{eqnarray}
\label{delta_0(q_2((x, r)), q_2((y, t))) = ...}
\lefteqn{\delta_0(q_2((x, r)), q_2((y, t)))} \nonumber \\
 & = & \inf\bigg\{\sum_{j = 1}^n
        \delta(q_2((x_j, r_j)), q_2((x_{j + 1}, r_{j + 1}))) : \nonumber \\
 & & \qquad\qquad  x_1, \ldots, x_{n + 1} \in X, \
                 r_1, \ldots, r_{n + 1} \in {\bf R} \hbox{ satisfy 
(\ref{q_2((x_1, r_1)) = q_2((x, r)), q_2(x_{n + 1}, r_{n + 1})) = q_2((y, t))})}                     \biggr\},
\end{eqnarray}
so that the infimum is taken over all finite sequences of elements of
$Y_2$ connecting $q_2((x, r))$ to $q_2((y, t))$.  In particular,
\begin{equation}
\label{delta_0(q_2((x, r)), q_2((y, t))) le delta(q_2((x, r)), q_2((y, t)))}
        \delta_0(q_2((x, r)), q_2((y, t))) \le \delta(q_2((x, r)), q_2((y, t))),
\end{equation}
since one can take $n = 2$, $x_1 = x$, $r_1 = r$, $x_2 = y$, and $r_2
= t$.  Of course, (\ref{delta_0(q_2((x, r)), q_2((y, t))) = ...}) is
nonnegative and symmetric in $q_2((x, r))$ and $q_2((y, t))$, because of
the corresponding properties of $\delta(\cdot, \cdot)$.  By construction,
$\delta_0(\cdot, \cdot)$ satisfies the triangle inequality
\begin{eqnarray}
\label{delta_0(q_2((x, r)), q_2((z, u))) le ...}
\lefteqn{\delta_0(q_2((x, r)), q_2((z, u)))} \\
 & \le & \delta_0(q_2((x, r)), q_2((y, t))) + \delta_0(q_2((y, t)), q_2((z, u)))
                                                                 \nonumber
\end{eqnarray}
for every $x, y, z \in X$ and $r, t, u \in {\bf R}$.  This is
basically because any finite sequences of elements of $Y_2$ connecting
$q_2((x, r))$ to $q_2((y, t))$ and connecting $q_2((y, t))$ to
$q_2((z, u))$ can be combined to get a finite sequence of elements of
$Y_2$ connecting $q_2((x, u))$ to $q_2((z, t))$.

        Let $q_2((x, r))$, $q_2((y, t))$ be any two elements of $Y_2$
again, and let (\ref{q_2((x_1, r_1)), ldots, q_2((x_{n + 1}, r_{n +
    1}))}) be a finite sequence of elements of $Y_2$ that satisfies
(\ref{q_2((x_1, r_1)) = q_2((x, r)), q_2(x_{n + 1}, r_{n + 1})) = q_2((y, t))}).
Observe that
\begin{equation}
\label{sum_{j = 1}^n delta(. , .) ge sum_{j = 1}^n dist(r_j - r_{j + 1}, Z)}
        \sum_{j = 1}^n \delta(q_2((x_j, r_j)), q_2((x_{j + 1}, r_{j + 1})))
          \ge \sum_{j = 1}^n \dist(r_j - r_{j + 1}, {\bf Z}),
\end{equation}
by (\ref{delta(q_2((x, r)), q_2((y, t))) ge dist(r - t, {bf Z})}).
The triangle inequality (\ref{dist(a + b, {bf Z}) le dist(a, {bf Z}) +
  dist(b, {bf Z})}) for $\dist(a, {\bf Z})$ implies that
\begin{equation}
\label{sum_{j = 1}^n dist(r_j - r_{j + 1}, Z) ge ... = dist(r - t, Z)}
 \sum_{j = 1}^n \dist(r_j - r_{j + 1}, {\bf Z}) \ge \dist(r_1 - r_{n + 1}, {\bf Z})
                                               = \dist(r - t, {\bf Z}),
\end{equation}
using the fact that $r_1 \equiv r$ and $r_{n + 1} \equiv t$ modulo
${\bf Z}$ in the second step, which follows from (\ref{q_2((x_1, r_1))
  = q_2((x, r)), q_2(x_{n + 1}, r_{n + 1})) = q_2((y, t))}).  Thus
\begin{equation}
\label{sum_{j = 1}^n delta(q_2(. , .) , q_2(. , .)) ge dist(r - t, {bf Z})}
 \sum_{j = 1}^n \delta(q_2((x_j, r_j)), q_2((x_{j + 1}, r_{j + 1})))
                 \ge \dist(r - t, {\bf Z}),
\end{equation}
and hence
\begin{equation}
\label{delta_0(q_2((x, r)), q_2((y, t))) ge dist(r - t, {bf Z})}
        \delta_0(q_2((x, r)), q_2((y, t))) \ge \dist(r - t, {\bf Z}),
\end{equation}
by taking the infimum of the sums on the left side of (\ref{sum_{j =
    1}^n delta(q_2(. , .) , q_2(. , .)) ge dist(r - t, {bf Z})}).

        As before, any two elements of $Y_2$ can be represented as
$q_2((x, r))$, $q_2((y, t))$ for some $x, y \in X$ and $r, t \in {\bf R}$
such that
\begin{equation}
\label{|r - t| le frac{1}{2} and |frac{r + t}{2}| le frac{1}{2}}
        |r - t| \le \frac{1}{2} \quad \hbox{and} \quad 
                               \biggl|\frac{r + t}{2}\biggr| \le \frac{1}{2},
\end{equation}
and hence $|r|, |t| \le 3/4$.  Let (\ref{q_2((x_1, r_1)), ldots,
  q_2((x_{n + 1}, r_{n + 1}))}) be a finite sequence of elements of
$Y_2$ such that (\ref{q_2((x_1, r_1)) = q_2((x, r)), q_2(x_{n + 1},
  r_{n + 1})) = q_2((y, t))}) again.  We may as well choose $x_1,
\ldots, x_{n + 1} \in X$ and $r_1, \ldots, r_{n + 1} \in {\bf R}$ such
that $x_1 = x$, $r_1 = r$, and
\begin{equation}
\label{|r_j - r_{j + 1}| = dist(r_j - r_{j + 1}, {bf Z})}
        |r_j - r_{j + 1}| = \dist(r_j - r_{j + 1}, {\bf Z})
\end{equation}
for $j = 1, \ldots, n$.  Suppose for the moment that
\begin{equation}
\label{sum_{j = 1}^n dist(r_j - r_{j + 1}, {bf Z}) < 1/2}
        \sum_{j = 1}^n \dist(r_j - r_{j + 1}, {\bf Z}) < 1/2,
\end{equation}
so that
\begin{equation}
\label{sum_{j = 1}^n |r_j - r_{j + 1}| < 1/2}
        \sum_{j = 1}^n |r_j - r_{j + 1}| < 1/2.
\end{equation}
In particular,
\begin{equation}
\label{|r - r_{n + 1}| = |r_1 - r_{n + 1}| < 1/2}
        |r - r_{n + 1}| = |r_1 - r_{n + 1}| < 1/2,
\end{equation}
which implies that
\begin{equation}
\label{|r_{n + 1} - t| le |r_{n + 1} - r_1| + |r - t| < 1/2 + 1/2 = 1}
        |r_{n + 1} - t| \le |r_{n + 1} - r_1| + |r - t| < 1/2 + 1/2 = 1,
\end{equation}
by the first part of (\ref{|r - t| le frac{1}{2} and |frac{r + t}{2}|
  le frac{1}{2}}).  It follows that $r_{n + 1} = t$ under these
conditions, since $r_{n + 1} \equiv t$ modulo ${\bf Z}$, by
(\ref{q_2((x_1, r_1)) = q_2((x, r)), q_2(x_{n + 1}, r_{n + 1})) =
  q_2((y, t))}).  Using (\ref{sum_{j = 1}^n dist(r_j - r_{j + 1}, {bf
    Z}) < 1/2}) again, we get that
\begin{equation}
\label{|r - r_l| + |r_l - t| < 1/2}
        |r - r_l| + |r_l - t| < 1/2
\end{equation}
for each $l = 1, \ldots, n + 1$, and hence
\begin{equation}
\label{|r_l - frac{r + t}{2}| le frac{|r_l - r| + |r_l - t|}{2} < frac{1}{4}}
 \biggl|r_l - \frac{r + t}{2}\biggr| \le \frac{|r_l - r| + |r_l - t|}{2} 
                                                               < \frac{1}{4}.
\end{equation}
Thus $|r_l| < 3/4$ for each $l = 1, \ldots, n + 1$, and of course
$|r_j - r_{j + 1}| < 1/2$ for each $j = 1, \ldots, n$, by (\ref{sum_{j
    = 1}^n |r_j - r_{j + 1}| < 1/2}).  This permits us to apply
(\ref{delta(q_2((x, r)), q_2((y, t))) ge C^{-1} rho((x, r), (y, t))})
to $q_2((x_j, r_j))$ and $q_2((x_{j + 1}, r_{j + 1}))$ for each $j =
1, \ldots, n$, to get that
\begin{equation}
\label{delta(q_2((x_j, r_j)), q_2((x_{j + 1}, r_{j + 1}))) ge ...}
        \delta(q_2((x_j, r_j)), q_2((x_{j + 1}, r_{j + 1}))) 
                        \ge C^{-1} \, \rho((x_j, r_j), (x_{j + 1}, r_{j + 1}))
\end{equation}
for each $j = 1, \ldots, n$.  Because $\rho(\cdot, \cdot)$ is a metric
on $X \times {\bf R}$, and hence satisfies the triangle inequality,
we get that
\begin{eqnarray}
\label{sum_{j = 1}^n delta(q_2(. , .), q_2(. , .)) ge ...}
\sum_{j = 1}^n \delta(q_2((x_j, r_j)), q_2((x_{j + 1}, r_{j + 1}))) & \ge & 
  C^{-1} \, \sum_{j = 1}^n \rho((x_j, r_j), (x_{j + 1}, r_{j + 1})) \nonumber \\
                                    & \ge & C^{-1} \, \rho((x, r), (y, t)).
\end{eqnarray}

        Otherwise, if (\ref{sum_{j = 1}^n dist(r_j - r_{j + 1}, {bf Z}) < 1/2})
does not hold, then
\begin{equation}
\label{sum_{j = 1}^n delta(q_2((x_j, r_j)), q_2((x_{j + 1}, r_{j + 1}))) ge 1/2}
 \sum_{j = 1}^n \delta(q_2((x_j, r_j)), q_2((x_{j + 1}, r_{j + 1}))) \ge 1/2,
\end{equation}
by (\ref{sum_{j = 1}^n delta(. , .) ge sum_{j = 1}^n dist(r_j - r_{j + 1}, Z)}).
Combining this with the previous case, and taking the infimum of the sums
on the left side, we get that
\begin{equation}
\label{delta_0(q_2((x, r)), q_2((y, t))) ge ...}
        \delta_0(q_2((x, r)), q_2((y, t)))
                  \ge \min\big(C^{-1} \, \rho((x, r), (y, t)), 1/2\big)
\end{equation}
when $r$ and $t$ satisfy (\ref{|r - t| le frac{1}{2} and |frac{r +
    t}{2}| le frac{1}{2}}).  Note that we also have
\begin{equation}
\label{delta_0(q_2((x, r)), q_2((y, t))) le rho((x, r), (y, t))}
        \delta_0(q_2((x, r)), q_2((y, t))) \le \rho((x, r), (y, t))
\end{equation}
under these conditions, by (\ref{delta(q_2((x, r)), q_2((y, t))) le
  rho((x, r), (y, t))}) and (\ref{delta_0(q_2((x, r)), q_2((y, t))) le
  delta(q_2((x, r)), q_2((y, t)))}).
 
        Suppose now that $X$ is bounded with diameter less than or equal 
to $k$ for some $k \ge 1/2$, so that
\begin{equation}
\label{d(x, y) le k, 2}
        d(x, y) \le k
\end{equation}
for every $x, y \in X$.  Of course, this can always be arranged by
replacing $d(x, y)$ with the minimum of $d(x, y)$ and $k$, as in the
previous section.  Alternatively, if $X$ is already bounded with
respect to $d(x, y)$, then one can get this condition by multiplying
$d(x, y)$ by a suitable positive constant.  In both cases, one can
check that the bilipschitz condition for $\phi$ would be maintained.

        Using (\ref{d(x, y) le k, 2}) and the definition (\ref{rho((x, r), 
(y, t)) = max(d(x, y), |r - t|), 2}) of $\rho(\cdot, \cdot)$, we get that
\begin{equation}
\label{rho((x, r), (y, t)) le k}
        \rho((x, r), (y, t)) \le k
\end{equation}
for every $x, y \in X$ and $r, t \in {\bf R}$ such that $|r - t| \le
1/2$.   This implies that
\begin{equation}
\label{delta_0(q_2((x, r)), ((y, t))) le delta(q_2((x, r)), q_2((y, t))) le k}
  \delta_0(q_2((x, r)), q_2((y, t))) \le \delta(q_2((x, r)), q_2((y, t))) \le k,
\end{equation}
for every $x, y \in X$ and $r, t \in {\bf R}$, because of
(\ref{delta_0(q_2((x, r)), q_2((y, t))) le delta(q_2((x, r)), q_2((y,
  t)))}) and the definition (\ref{delta(q_2((x, r)), q_2((y, t))) =
  ...}) of $\delta(q_2((x, r)), q_2((y, t)))$.  If $r$, $t$ satisfy
(\ref{|r - t| le frac{1}{2} and |frac{r + t}{2}| le frac{1}{2}}), then
we get that
\begin{equation}
\label{rho((x, r), (y, t)) le max(C, 2 k) delta_0((q_2((x, r)), q_2((y, t)))}
        \rho((x, r), (y, t))
               \le \max(C, 2 \, k) \, \delta_0((q_2((x, r)), q_2((y, t))),
\end{equation}
by combining (\ref{delta_0(q_2((x, r)), q_2((y, t))) ge ...}) and
(\ref{rho((x, r), (y, t)) le k}).

\section{Nonnegative Borel measures}
\label{nonnegative borel measures}
\setcounter{equation}{0}

        Let us return to the setting of Section \ref{a basic situation}.  
If $\mu$ is a nonnegative Borel measure on $X$, then we would like to
have a corresponding product measure on $X \times {\bf R}$, using
Lebesgue measure on ${\bf R}$.  Of course, the standard product
measure construction applies when $\mu$ is at least $\sigma$-finite on
$X$.  It is better for $X$ to also have a countable base for its
topology, so that there is a countable base for the topology of $X
\times {\bf R}$ consisting of products of open subsets of $X$ and
${\bf R}$.  This implies that open subsets of $X \times {\bf R}$ are
measurable with respect to the usual product $\sigma$-algebra, and
hence that Borel sets in $X \times {\bf R}$ are measurable too.
Alternatively, if $X$ is a locally compact Hausdorff space, then one
might start with a Borel measure $\mu$ on $X$ with suitable regularity
properties, and look for a product Borel measure on $X \times {\bf R}$
with similar regularity properties.  More precisely, one can view
this in terms of nonnegative linear functionals on spaces of
continuous functions with compact support, using the Riesz
representation theorem.

        At any rate, such a product measure on $X \times {\bf R}$
is invariant under translations on ${\bf R}$, because Lebesgue measure
is invariant under translations.  If $\mu$ is invariant under $\phi$
on $X$, then the product measure is invariant under $\Phi$ on 
$X \times {\bf R}$.  One can then localize to get a measure on $Y_2$
that is invariant under the mappings $A_r$ corresponding to translation
on ${\bf R}$.  If $\mu$ is not invariant under $\phi$, then one can
still get measures on $Y_2$, by restricting the product measure to
the product of $X$ with an interval $I$ in ${\bf R}$ with length $1$.
Of course, the resulting measures on $Y_2$ will depend on $I$, but under
suitable conditions on $\phi$ and $\mu$, they may be reasonably similar.

        Suppose that $X$ is compact, and that the topology on $X$ is
determined by a metric $d(x, y)$.  If $\phi$ is bilipschitz with
respect to this metric, then one can get a metric on $Y_2$ that looks
locally approximately like a product metric on $X \times {\bf R}$, as
in the previous section.  If $\mu$ is Ahlfors regular of some
dimension $t$, then one can get an Ahlfors regular measure on $Y_2$ of
dimension $t + 1$, even if $\mu$ is not invariant under $\phi$.  More
precisely, if $\mu$ is Ahlfors regular on $X$ of dimension $t$, then
$\mu$ is approximately the same as $t$-dimensional Hausdorff measure
on $X$, in the sense that each is bounded by constant multiples of the
other.  Hausdorff measure in any dimension is approximately preserved
to within bounded factors by a bilipschitz mapping, which implies that
$\mu$ is approximately preserved by $\phi$ to within bounded factors
under these conditions.

        Even if $\mu$ is not Ahlfors regular, it may be approximately
preserved to within bounded factors by $\phi$, so that one can get
corresponding measures on $Y_2$ that are at least comparable to each
other.  If $\mu$ is a doubling measure on $X$, and if $\phi$ is
bilipschitz or at least quasisymmetric on $X$, then $\mu$ is
transformed by $\phi$ to a doubling measure on $X$, but the new
measure may not be comparable to $\mu$.  If $\mu$ is a doubling measure
on $X$ which is approximately preserved to within bounded factors by
$\phi$, and if $\phi$ is bilipschitz, then one can get doubling measures
on $Y_2$ from $\mu$, as before.  Although quasisymmetry of $\phi$ on $X$
is a very natural geometric condition, it is not by itself so convenient
for looking at geometric structures on $X \times {\bf R}$, and hence $Y_2$.
However, if $\phi$ is a quasisymmetric mapping on $X$ that approximately
preserves a nontrivial doubling measure $\mu$ on $X$ to within bounded
factors, then one can use that to get another geometric structure on $X$
that is approximately preserved by $\phi$ to within bounded factors,
at least if $X$ is also uniformly perfect.

\section{Some examples}
\label{some examples}
\setcounter{equation}{0}

        Let $B$ be a finite set with at least two elements, and let $X$
be the set of doubly-infinite sequences $x = \{x_j\}_{j = -\infty}^\infty$
such that $x_j \in B$ for each $j \in {\bf Z}$.  Equivalently, $X$
is the Cartesian product of a family of copies of $B$, indexed by
${\bf Z}$.  Thus $X$ is a compact Hausdorff topological space, with
respect to the product topology associated to the discrete topology
on each copy of $B$.  Let $\phi$ be the shift mapping defined by
\begin{equation}
\label{phi(x) = {x_{j - 1}}_{j = -infty}^infty}
        \phi(x) = \{x_{j - 1}\}_{j = -\infty}^\infty,
\end{equation}
which is a homeomorphism from $X$ onto itself.  Also let $\Phi$,
$Y_2$, etc., be as in Section \ref{a basic situation}, using this $X$
and $\phi$.

        Suppose that $E_0$ is a nonempty open set in $X$, and let $x$
be an element of $E_0$.  Because of the way that the product topology
is defined on $X$, there is a nonnegative integer $n$ such that $E_0$
contains every $y \in X$ that satisfies
\begin{equation}
\label{y_j = x_j for each j in {bf Z} with |j| le n}
        y_j = x_j \hbox{ for each } j \in {\bf Z} \hbox{ with } |j| \le n.
\end{equation}
If $\phi(E_0) = E_0$, then $\phi^k(E_0) = E_0$ for every $k \in {\bf
  Z}$, and hence $E_0$ also contains every $z \in X$ such that $y =
\phi^k(z)$ satisfies (\ref{y_j = x_j for each j in {bf Z} with |j| le
  n}) for some $k \in {\bf Z}$.  If $E_0$ is a closed set in $X$ too,
then it follows that $E_0 = X$, again because of the way that the
product topology is defined on $X$.  This shows that $E_0 = X$ when
$E_0$ is a nonempty open and closed subset of $X$ that is invariant
under $\phi$, which implies that $Y_2$ is connected, as in Section
\ref{connectedness}.

        Let $w$ be a nonnegative real-valued function on $B$ such that
\begin{equation}
\label{sum_{b in B} w(b) = 1}
        \sum_{b \in B} w(b) = 1.
\end{equation}
Thus $w$ defines a probability measure on $B$ in the obvious way, and
we let $\mu = \mu_w$ be the corresponding product measure on $X$,
using the probability measure on $B$ associated to $w$ on each factor.
One can first define the corresponding integral of a continuous
real-valued function on $X$ as a limit of suitable finite sums, and
then get $\mu_w$ as a regular Borel measure on $X$ from the Riesz
representation theorem.  Of course, $\mu_w$ is invariant under $\phi$
for every $w$, since $\mu_w$ is defined using the same probability
measure on each copy of $B$ in the product.  Note that there is a
countable base for the topology of $X$, by standard arguments.

        Let $x, y \in X$ be given, with $x \ne y$, and let $n(x, y)$
be the largest nonnegative integer such that
\begin{equation}
\label{x_j = y_j for every j in {bf Z} with -n + 1 le j le n}
 x_j = y_j \hbox{ for every } j \in {\bf Z} \hbox{ with } -n + 1 \le j \le n,
\end{equation}
which holds trivially when $n = 0$.  If $x = y$, then one can put
$n(x, y) = +\infty$.  It is easy to see that
\begin{equation}
\label{n(x, y) = n(y, x)}
        n(x, y) = n(y, x)
\end{equation}
for every $x, y \in X$, and that
\begin{equation}
\label{n(x, z) ge min(n(x, y), n(y, z))}
        n(x, z) \ge \min(n(x, y), n(y, z))
\end{equation}
for every $x, y, z \in X$.  Observe also that
\begin{equation}
\label{n(x, y) - 1 le n(phi(x), phi(y)) le n(x, y) + 1}
        n(x, y) - 1 \le n(\phi(x), \phi(y)) \le n(x, y) + 1
\end{equation}
for every $x, y \in X$.

        Let $a$ be a positive real number which is strictly less than $1$,
and put
\begin{equation}
\label{d_a(x, y) = a^{n(x, y)}}
        d_a(x, y) = a^{n(x, y)}
\end{equation}
for every $x, y \in X$, which is interpreted as being equal to $0$ when
$x = y$.  Thus $d_a(x, y) > 0$ when $x \ne y$,
\begin{equation}
\label{d_a(x, y) = d_a(y, x)}
        d_a(x, y) = d_a(y, x)
\end{equation}
for every $x, y \in X$, and
\begin{equation}
\label{d_a(x, z) le max(d_a(x, y), d_a(y, z))}
        d_a(x, z) \le \max(d_a(x, y), d_a(y, z))
\end{equation}
for every $x, y, z \in X$, by (\ref{n(x, y) = n(y, x)}) and (\ref{n(x,
  z) ge min(n(x, y), n(y, z))}).  This shows that $d_a(x, y)$ is a
metric on $X$ for each $a \in (0, 1)$, and in fact $d_a(x, y)$ is an
ultrametric on $X$, since it satisfies the ultrametric version
(\ref{d_a(x, z) le max(d_a(x, y), d_a(y, z))}) of the triangle
inequality.  By construction, the topology on $X$ determined by
$d_a(x, y)$ is the same as the product topology described earlier for
each $a \in (0, 1)$.  In particular, these ultrametrics on $X$ are
topologically equivalent, and indeed we have that
\begin{equation}
\label{d_a(x, y)^alpha = d_{a^alpha}(x, y)}
        d_a(x, y)^\alpha = d_{a^\alpha}(x, y)
\end{equation}
for every $a \in (0, 1)$, $\alpha > 0$, and $x, y \in X$.

        It follows from (\ref{n(x, y) - 1 le n(phi(x), phi(y)) le n(x, y) + 1})
that
\begin{equation}
\label{a d_a(x, y) le d_a(phi(x), phi(y)) le (1/a) d_a(x, y)}
        a \, d_a(x, y) \le d_a(\phi(x), \phi(y)) \le (1/a) \, d_a(x, y)
\end{equation}
for every $a \in (0, 1)$ and $x, y \in X$, so that $\phi$ is
bilipschitz with constant $C = 1/a$ with respect to $d_a(x, y)$.
However, one can also check that the collection of iterates $\phi^k$
of $\phi$ with $k \in {\bf Z}$ is not equicontinuous at any point in
$X$ with respect to $d_a(x, y)$ for any $a \in (0, 1)$.  If $d(x, y)$
is any metric on $X$ that determines the same topology on $X$, then
the identity mapping on $X$ is uniformly continuous as a mapping from
$X$ equipped with $d(x, y)$ into $X$ equipped with $d_a(x, y)$ for any
$a \in (0, 1)$, because $X$ is compact.  This implies that the collection
of iterates $\phi^k$ of $\phi$ with $k \in {\bf Z}$ is not equicontinuous
with respect to any metric $d(x, y)$ on $X$ that determines the same
topology on $X$.

        By construction,
\begin{equation}
\label{d_a(x, y) le 1}
        d_a(x, y) \le 1
\end{equation}
for each $x, y \in X$ and $a \in (0, 1)$, and equality holds when $x_0
\ne y_0$.  The closed ball in $X$ with respect to $d_a(x, y)$ centered
at some point $x \in X$ and with radius $a^n$ for some nonnegative
integer $n$ is the same as the set of $y \in X$ that satisfy (\ref{x_j
  = y_j for every j in {bf Z} with -n + 1 le j le n}).  If $w$ is a
positive real-valued function on $B$ that satisfies (\ref{sum_{b in B}
  w(b) = 1}), then one can check that the corresponding probability
measure $\mu_w$ on $X$ is a doubling measure with respect to $d_a(x,
y)$.

        Suppose now that $w$ corresponds to the uniform distribution on
$B$, so that
\begin{equation}
\label{w(b) = 1/num(B)}
        w(b) = 1/\# B,
\end{equation}
where $\# B$ denotes the number of elements of $B$.  In this case, the
measure with respect to $\mu_w$ of a closed ball in $X$ with respect
to $d_a(x, y)$ with radius $a^n$ for some nonnegative integer $n$ is
\begin{equation}
\label{(num(B))^{-2 n}}
        (\# B)^{-2 n}.
\end{equation}
If we put $t = - 2 \log (\# B) / \log a$, then $t > 0$ and
\begin{equation}
\label{(a^n)^t = (num(B))^{-2 n}}
        (a^n)^t = (\# B)^{-2 n}
\end{equation}
for each $n \ge 0$, which implies that $\mu_w$ is Ahlfors regular on
$X$ with respect to $d_a(x, y)$, with dimension $t$.  In particular,
the Hausdorff dimension of $X$ with respect to $d_a(x, y)$ is equal to
$t$.

\section{Snowflake metrics}
\label{snowflake metrics}
\setcounter{equation}{0}

        If $d(x, y)$ is a metric on a set $X$, then $d(x, y)^\alpha$
is also a metric on $X$ when $0 < \alpha < 1$, which determines the
same topology on $X$.  However, $d(x, y)^\alpha$ does not necessarily
satisfy the triangle inequality when $\alpha > 1$, even when $X$ is
the unit interval with the standard metric.  It is easy to see that
$d(x, y)^\alpha$ is still a quasi-metric on $X$ when $\alpha > 1$,
which means that it satisfies a weaker version of the triangle
inequality with an extra constant factor on the right side, and which
is adequate in many situations.  Of course, if $d(x, y)$ is an
ultrametric on $X$, then $d(x, y)^\alpha$ is an ultrametric on $X$ for
each $\alpha > 0$.  Note that the Hausdorff dimension of $X$ with
respect to $d(x, y)^\alpha$ is equal to the Hausdorff dimension of
$X$ with respect to $d(x, y)$ divided by $\alpha$.

        If $\phi$ is a bilipschitz mapping from $X$ onto itself with
respect to $d(x, y)$ with constant $C$, then $\phi$ is also
bilipschitz with respect to $d(x, y)^\alpha$, with constant
$C^\alpha$.  In particular, if $\phi$ is an isometry with respect to
$d(x, y)$, then $\phi$ preserves $d(x, y)^\alpha$ for each $\alpha$.
Thus one can repeat the same types of constructions as before, with
$d(x, y)$ replaced with $d(x, y)^\alpha$.  This was already built in
the examples discussed in the previous section, using the parameter $a
\in (0, 1)$.  If $X = [0, 1]$ with the standard metric $d(x, y) = |x - y|$
and $0 < \alpha < 1$, then $X$ is basically a snowflake curve of dimension 
$1/\alpha$ with respect to $d(x, y)^\alpha$.

        If $d(x, y)$ is any metric on a set $X$ and $0 < \alpha < 1$,
then one can check that every continuous curve in $X$ with finite
length with respect to $d(x, y)^\alpha$ is constant.  Consider the
metric on $X \times {\bf R}$ defined by
\begin{equation}
\label{rho_alpha((x, r), (y, t)) = max(d(x, y)^alpha, |r - t|)}
        \rho_\alpha((x, r), (y, t)) = \max(d(x, y)^\alpha, |r - t|),
\end{equation}
which is the analogue of (\ref{rho((x, r), (y, t)) = max(d(x, y), |r -
  t|)}), (\ref{rho((x, r), (y, t)) = max(d(x, y), |r - t|), 2}) with
$d(x, y)$ replaced by $d(x, y)^\alpha$.  If $\gamma$ is any continuous
curve in $X \times {\bf R}$ with finite length with respect to
(\ref{rho_alpha((x, r), (y, t)) = max(d(x, y)^alpha, |r - t|)}), then
the projection of $\gamma$ into $X$ has finite length with respect to
$d(x, y)^\alpha$, and hence is constant.  This means that $\gamma$ is
contained in a line $\{x\} \times {\bf R}$ for some $x \in X$, and
hence corresponds to a curve of finite length in the real line, with
the standard metric.  Similarly, if $Y_2$ is equipped with a metric
that looks locally like (\ref{rho_alpha((x, r), (y, t)) = max(d(x,
  y)^alpha, |r - t|)}), as before, then any continuous curve in $Y_2$
with finite length has to be contained in the image of a line $\{x\}
\times {\bf R}$ under the usual quotient mapping $q_2$.

\end{document}